\documentclass[final, 3p, times]{elsarticle}

\usepackage{hyperref, float, amsmath, bm, listings, color, footnote}
\usepackage{amsfonts, amssymb, comment, graphicx, subcaption}
\usepackage[linesnumbered,ruled,vlined]{algorithm2e}
\usepackage[bottom]{footmisc}
\usepackage[dvipsnames]{xcolor}

\newcommand{\lexint}{\texttt{LeXInt}\xspace}

\hypersetup{colorlinks,
            linkcolor = {red!50!black},
            citecolor = {red!10!black},
            urlcolor  = {green!10!black}}

\definecolor{codegreen}{rgb}{0.1, 0.1, 0.9}
\definecolor{codegray}{rgb}{0.9, 0, 0.5}
\definecolor{codepurple}{rgb}{0.85, 0.2, 0.2}
\definecolor{backcolour}{rgb}{0.97, 0.97, 0.97}

\lstdefinestyle{CodeColour}{
	backgroundcolor=\color{backcolour},   
	commentstyle=\color{codegreen},
	keywordstyle=\color{codepurple},
	numberstyle=\tiny\color{codegray},
	stringstyle=\color{codepurple},
	basicstyle=\ttfamily\footnotesize,
	breakatwhitespace=false,         
	breaklines=true,                 
	captionpos=b,                    
	keepspaces=true,                 
	numbers=left,                    
	numbersep=5pt,                  
	showspaces=false,                
	showstringspaces=false,
	showtabs=false,                  
	tabsize=4
}

\journal{SoftwareX}

\begin{document}

\begin{frontmatter}



\title{\lexint: Package for Exponential Integrators employing Leja interpolation}

\author[inst1]{Pranab J. Deka\corref{lod1}} \ead{pranab.deka@uibk.ac.at} 
\author[inst1]{Lukas Einkemmer}             \ead{lukas.einkemmer@uibk.ac.at}
\author[inst2]{Mayya Tokman}                \ead{mtokman@ucmerced.edu}

\cortext[lod1]{Corresponding author}

\affiliation[inst1]{organization = {Department of Mathematics, University of Innsbruck},
            city = {Innsbruck},
            postcode = {6020}, 
            country = {Austria}}
            
\affiliation[inst2]{organization = {School of Natural Sciences, University of California},
            city = {Merced},
            postcode = {CA 95343}, 
            country = {USA}}

\begin{abstract}
We present a publicly available software for exponential integrators that computes the $\varphi_l(z)$ functions using polynomial interpolation. The interpolation method at Leja points have recently been shown to be competitive with the traditionally-used Krylov subspace method. The developed framework facilitates easy adaptation into any \textsc{Python} software package for time integration. 
\end{abstract}


% \begin{highlights}
% \end{highlights}


\end{frontmatter}


%%%%%%%%%%%%%%%%%%%%%%%%%%%%%%%%%%%%%%%%%%%%%%%%%%%%%%%%%%%%%%%%%%%%%%%%%%%%%%%%%%


\section{Motivation and significance}

Time-dependent partial differential equations (PDEs) are ubiquitous in various fields in science. Integrating PDEs in time with high accuracy whilst incurring as little computational cost as possible is highly desirable. Substantial amount of research has been devoted to the development of numerical algorithms and codes to perform high-resolution simulations with high fidelity. 

Explicit temporal integrators are widely used in many scenarios owing to the simplicity of their algorithm and implementation. However, as the number of physical processes considered in a certain PDE increases or if the stiff nature of the underlying PDE becomes prominent, the performance of the explicit integrators is severely deteriorated owing to the stability constraints. The increase in the stiffness of an equation results in ever more stringent Courant--Friedrich--Levy (CFL) time step size limit. Implicit integrators have been widely used as alternatives to the explicit methods owing to their ability to take large step sizes. They can provide substantial boost to the simulations. In many practical cases, however, one has to resort to iterative schemes to solve large systems of linear equations. Furthermore, the use of preconditioners to speed-up the simulations is a common practice in many situations. In some cases, the simulations fail to converge without the use of a good preconditioner. The complexity involved in such an algorithm may make them unfavourable for intricate problems.

Exponential integrators are a class of temporal integrators that linearise the underlying PDE at every time step - the linear term is solved exactly (in time) and the nonlinear term is approximated with some explicit methods. An extensive review on exponential integrators have been presented by Hockbruck \& Ostermann \cite{Ostermann10}. Let us consider the autonomous system
\begin{equation*}
    \frac{\partial u}{\partial t} = f(u(t)), \qquad u(t = 0) = u^0,
\end{equation*}
where $f(u)$ is some nonlinear function of $u$. We re-write the above equation as
\begin{equation*}
	\frac{\partial u}{\partial t} = \mathcal{A} \, u + g(u),
\end{equation*}
where $\mathcal{A}$ is a matrix (usually the stiff linear part of $f(u)$) and $g(u)$ is the nonlinear part. An approximation to the solution of this equation is given by
\begin{equation*}
    u^{n + 1} = u^n + \varphi_1(\mathcal{A} \Delta t) f(u^n) \Delta t,
\end{equation*}
where $u^n$ is the solution at the $n^\mathrm{th}$ time step. This is the first-order exponential Euler method. If $\mathcal{A}$ is replaced by the Jacobian of $f(u)$ evaluated at the respective time step, one obtains the Rosenbrock--Euler method, which is of second order. It is to be noted that replacing the linear part by the Jacobian allows for one to obtain higher-order schemes with fewer $\varphi_l(z)$ function evaluations. The $\varphi_l(z)$ functions are given by:
\begin{equation*}
    \varphi_{l + 1}(z) = \frac{1}{z} \left(\varphi_l(z) - \frac{1}{l!} \right),  \quad l \geq 0,
\end{equation*}
where $\varphi_0(z) = e^z$ corresponds to the matrix exponential. We compute these $\varphi_l(z)$ functions,  the most expensive part of exponential integrators, using the method of polynomial interpolation at Leja points. Details on this iterative scheme are provided in Sec. \ref{sec:software}.

Exponential integrators do not suffer from any CFL restrictions (unlike explicit integrators), are unconditionally stable, and can take much larger step sizes than implicit methods. This makes them highly attractive for solving time-dependent problems. Additionally, one can obtain the exact solution of a linear autonomous PDE (subject to the spatial discretisation) for any given step size. This is an added bonus over the implicit integrators, as they will always incur some error, irrespective of their order of convergence.

We make our tools and contributions available to the scientific community with the release of the \textbf{Le}ja interpolation for e\textbf{X}ponential \textbf{Int}egrators (\lexint; \url{https://github.com/Pranab-JD/LeXInt}) package. This is a cumulation of the algorithms implemented and tested out in our previous work where we studied the performance of an automatic step-size controller for improved computational efficiency \cite{Deka22a} and analysed the performance of the Leja method with explicit, implicit, and Krylov-based exponential integrators for the set of equations of magnetohydrodynamics (MHD) \cite{Deka22b}. We provide a user-friendly framework (i.e. \textsc{Python}) for a range of exponential integrators as well as the Leja interpolation method. This is to introduce the scientific community to integrators that are well-suited for temporal integration of stiff as well as general-purpose problems and the highly effective iterative techniques that are used in such integrators. To the best of our knowledge, this is the first open-access Leja-method-based compilation of exponential integrators package.

%%%%%%%%%%%%%%%%%%%%%%%%%%%%%%%%%%%%%%%%%%%%%%%%%%%%%%%%%%%%%%%%%%%%%%%%%%%%%%%%%%

\section{Software description}
\label{sec:software}

\lexint comprises of several exponential integrators suited for both constant and variable step size implementation. The integrators are implemented in a modular format, in essence, any integrator can easily be integrated into the package and any integrator can be used for any given problem. However, it is to be noted that the performance of the integrators may vary with the problem under consideration. We primarily focus on integrators that are based on linearising the underlying PDE by computing the Jacobian at every time step: Exponential Rosenbrock (EXPRB) and Exponential Propagation Iterative Runge--Kutta (EPIRK) integrators. We have adopted the vertical implementation procedure, proposed by Rainwater \& Tokman \cite{Tokman16}, for optimised performance. Although this was initially proposed only for the Krylov subspace algorithm, the vertical approach can, very well, result in substantial amount of computational savings for the Leja interpolation method.

For adaptive step size implementation, one requires an error estimate at every time step. One of the cheapest ways to compute the error estimate is if it is inherently embedded in the integrator, i.e. an embedded integrator. This is why we focus only on embedded exponential integrators, where the error estimate does not require additional internal stages. The list of embedded exponential integrators implemented in \lexint include EXPRB32 \cite{Caliari09, Hochbruck09, Ostermann10}, EXPRB43 \cite{Caliari09, Hochbruck09, Ostermann10}, EXPRB53s3 \cite{Luan14}, EXPRB54s4 \cite{Luan14}, EPIRK4s3 \cite{Tokman17a, Tokman17b}, EPIRK4s3A \cite{Tokman16}, and EPIRK5P1 \cite{Tokman12}. Each of these integrators have been implemented in \lexint in a way that the integrator function returns the lower-order and the higher-order solutions. The difference between these two solutions yields an estimate of the error incurred. As there are a multitude of step-size controllers in the literature, we give the user full flexibility in choosing their desired step-sizing strategy. The function also returns the number of matrix-vector products computed at a given time step that can be considered as a proxy of the computational cost. Obviously, the cost of the matrix-vector products will depend on the size of the matrix and the vector, i.e. the number of grid points. However, we are only interested in the determining the cost of the algorithm and not on the dependence of the cost on the grid size. As such, we assume that the computational cost is normalised to the size to the matrix (and the vector). For constant step sizes, in addition to the aforementioned ones, we have implemented Rosenbrock--Euler \cite{Pope63}, EXPRB42 \cite{Luan17}, EPIRK4s3B \cite{Tokman16}, and EPIRK5P2 \cite{Tokman12}. In cases where the integrators do not possess an embedded error estimator, one can generate an error estimate using Richardson extrapolation.

The inputs to an integrator include the state variable(s), the right-hand-side (RHS) function, the step size, the scaling and the shifting factors (see following subsection), the set of predetermined Leja points, and the option to interpolate on real or imaginary Leja points depending on the problem under consideration. The interpolation method consists of computing the divided differences or the coefficients of the  polynomial to be interpolated. This depends on the step size, shifting and scaling factors, the integrator coefficients, and the order of the $\varphi_l(z)$ function. Once the coefficients are determined, the first term of the polynomial is computed by multiplying the function to be interpolated with the first coefficient. The subsequent terms are added to this polynomial one after another until the desired accuracy is reached (see following subsection).

%%% --------------------------------------------------------- %%%

\subsection*{Polynomial Interpolation at Leja points}

One of the crucial aspects of efficient implementation of exponential integrators is an adept iterative scheme. Whilst the Krylov subspace algorithm has long been proposed as an effective iterative scheme for exponential integrators \cite{Sidje98, Moler03, Ostermann10, Tokman12, Einkemmer17}, it does have the drawback of having the need to compute inner products that becomes a serious impediment on massively parallel structures (GPUs). Recently introduced adaptive Krylov-based \texttt{KIOPS} \cite{Gaudreault18} utilises incomplete orthogonalisation procedure which can reduce the number of inner products needed per Krylov iteration to as few as two. Thus, the use of incomplete orthogonalisation significantly alleviates the challenges associated with parallelisation of the computation of inner products. However, for problems where more Krylov vectors, and consequently, more inner products need to be computed, this issue persists. Additionally, even in the case of incomplete orthogonalisation, one needs to store a Krylov space which may be problematic for very large-sized problems.

We choose the method of polynomial interpolation at Leja points \cite{Leja1957, Reichel90, Baglama98} that has been shown to be competitive with the Krylov-based methods \cite{Bergamaschi06, Caliari07b, Deka22b}. This can be attributed mainly to the simplicity of the algorithm. One of the minor drawbacks of the Leja interpolation method is that it needs some approximation of the spectrum. It is to be noted that one needs only a crude estimate of the largest and the smallest eigenvalue of the matrix (for linear equations) or the Jacobian (for nonlinear equations). Explicitly forming the Jacobian matrix at every time step is computationally unattractive, which is why we only consider the action of the Jacobian matrix on the relevant vector. In the case of Jacobian-free (or even matrix-free) implementation, the method of power iterations can be employed to compute the spectrum every `n' time steps. In Deka \& Einkemmer \cite{Deka22b}, we have shown this to be an efficient technique for the highly nonlinear MHD equations. We note that the method of power iterations gives us only the magnitude of the largest eigenvalue. The nature of the problem, as in whether it is diffusion-dominated or advection-dominated, can be understood by comparing the diffusion and the advection CFL times. For diffusion-dominated problems, we assume that the largest eigenvalue (in magnitude) lies on the real axis and we set the smallest eigenvalue to 0, whereas for advection-dominated problems, we assume that the eigenvalues predominantly lie on the imaginary axis, and we interpolate the polynomial on imaginary Leja points. Here, the smallest eigenvalue is chosen to be the negative of the largest eigenvalue obtained using power iterations.

Now that we have an estimate of the largest ($\alpha$) and the smallest ($\beta$) eigenvalue (in magnitude), the scaling and the shifting factors can be defined as $c = (\alpha + \beta)/2$ and $\gamma = (\beta - \alpha)/4$, respectively \cite{Caliari04, Caliari14}. For real eigen values, $\beta = 0$, whereas for imaginary eigenvalues, $\beta = -\alpha$. The factor of $4$ emerges from the fact that we have chosen Leja points ($\xi$) in the arbitrary spectral domain $[-2, 2]$. We compute the coefficients ($d$) of the polynomial, to be interpolated, using the divided differences algorithm. Then, we form the polynomial by adding an additional term at every iteration until the desired accuracy is reached. This can be mathematically written as 
\begin{align*}
    p_{m+1}(z) & = p_m(z) + d_{m+1} \, y_{m+1}(z), \\
    y_{m+1}(z) & = y_m(z) \times \left(\frac{z - c}{\gamma} - \xi_{m} \right),
\end{align*}
where $d_m$ is the $m^\mathrm{th}$ coefficient of the polynomial, $\xi_m$ is the $m^\mathrm{th}$ Leja point, and $p_m(z)$ is the $m^\mathrm{th}$ term of the polynomial for $m \geq 0$. The polynomial is initialised as $p_0 = d_0\,y_0$. Let us clearly state that $p_m$ corresponds to the polynomial interpolation of $\exp(\Delta t(c + \gamma \xi))\,v$ if the divided differences have been computed for $\exp(\Delta t(c + \gamma \xi))$ whereas if the coefficients have been computed for $\varphi_l(\Delta t(c + \gamma \xi))$, $p_m$ would be the polynomial interpolation of $\varphi_l(\Delta t(c + \gamma \xi))\,v$ (for some vector $v$). If the values of the user-defined tolerances are too large or if the error incurred at a certain time step becomes minuscule, a step size controller may allow for extremely large step sizes. Such large step sizes might cause the algorithm to diverge, which is why, we adopt a safety measure of checking that the error incurred does not exceed a certain threshold. Very large value of this error is an indication of impending divergence of the algorithm. In such a case, we reject the step size and restart the time step with a smaller step size. We note that $z$ is the Jacobian-vector product of the underlying PDE evaluated at $u^n$ and it changes at every iteration depending on the value of $y_m(z)$. The convergence is determined by the error incurred (i.e. $|d_m| \, \|y_m\|$) being less than the user-defined tolerance (subject to a safety factor). 

\lexint has two functions for interpolating $\varphi_l(z)$ on real (\texttt{`real\_Leja\_phi'}) and imaginary (\texttt{`imag\_Leja\_phi'}) Leja points. To speed-up convergence, it is recommended that if the largest eigenvalue of the Jacobian under consideration is real, one interpolates the exponential-like function on real Leja points, whereas if the largest eigenvalue (in magnitude) is imaginary, the interpolation is performed on imaginary Leja points. If the magnitude of the largest real and imaginary eigenvalues are relatively similar, one could interpolate on either real or imaginary Leja points.

The exact solution (in time) of an autonomous homogeneous linear differential equation of the form \[ \frac{\partial u}{\partial t} = \mathcal{A}u, \] is given by $u^{n+1} = \exp({\mathcal{A} \Delta t) u^n}$. \lexint provides the functions \texttt{`real\_Leja\_exp'} and \texttt{`imag\_Leja\_exp'} to compute the matrix exponential by means of Leja interpolation on real and imaginary Leja points, respectively. In the case of autonomous nonhomogenous linear differential equations of the form \[ \frac{\partial u}{\partial t} = \mathcal{A}u + \mathcal{S}, \] where $\mathcal{S}$ is some forcing term, the exact solution, in time, is given by $u^{n+1} = u^n + \varphi_1(\mathcal{A} \Delta t) (\mathcal{A} u^n + S) \Delta t$. We have provided the functions \texttt{`real\_Leja\_phi\_nl'} or \texttt{`imag\_Leja\_phi\_nl'} to compute the (exact) solutions of such equations. One can compute the exact solution (in time) for autonomous linear equations, which is why the functions for interpolation of the matrix exponential are provided only for constant step size implementation, i.e., without any error estimate. The desired accuracy can be chosen by the user by tuning the tolerance. 

It is to be noted that \lexint can work in fully matrix-free structure as well as in any given formulation of the matrix, provided that one has a well-defined RHS function (this is similar to how one would implement an explicit method). Obviously, a matrix-free formulation is preferable from the computational viewpoint.

%%%%%%%%%%%%%%%%%%%%%%%%%%%%%%%%%%%%%%%%%%%%%%%%%%%%%%%%%%%%%%%%%%%%%%%%%%%%%%%%%%

\section{Illustrative examples}

\begin{figure}[tb]
    \centering
    \includegraphics[width = 0.7\columnwidth]{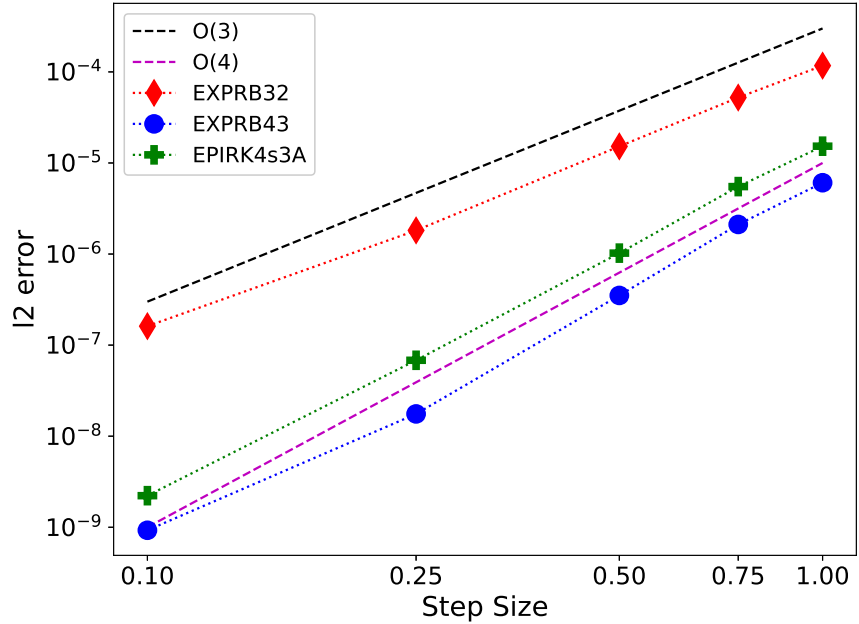}
    \caption{Convergence curves of EXPRB32, EXPRB43, and EPIRK4s3A.}
    \label{fig:order}
\end{figure}

We show the performance of a selected number of integrators for a couple of problems. These problems have been drawn from Refs. \cite{Einkemmer18, Deka22a}. In both these problems, we consider periodic boundary conditions on $[0, 1]$. We have used the second-order centered difference scheme to discretise the 1D Laplacian operator ($\partial^2/\partial x^2$) and third-order upwind scheme to discretise the advection term ($\partial/\partial x$). The first example is the Burgers' equation,
\begin{equation}
    \frac{\partial u}{\partial t} = \frac{\partial^2 u}{\partial x^2} + \frac{\eta}{2} \frac{\partial u^2}{\partial x}, 
    \label{eq:burgers}
\end{equation}
where $\eta$ is the P\'eclet number and the initial condition is given by
\begin{equation*}
	u(x, t = 0) = 1 + \exp\left(1 - \frac{1}{1-(2x - 1)^2}\right) + \frac{1}{2} \exp\left(-\frac{(x - x_0)^2}{2\sigma^2}\right),
\end{equation*}
with $x_0 = 0.9$ and $\sigma = 0.02$. We consider two different cases of the resolution, in terms of the number of grid points ($N$), $\eta$, and the simulation time $t_f$: (a) $N = 64, \eta = 200, t_f = 10^{-3}$ and (b) $N = 256, \eta = 10, t_f = 10^{-2}$. The second example is the Allen---Cahn equation:
\begin{equation}
    \frac{\partial u}{\partial t} = \frac{\partial^2 u}{\partial x^2} + 100 \, \left(u - u^3\right).
    \label{eq:ac}
\end{equation}
The initial condition is chosen to be
\begin{equation*}
    u(x, t = 0) = A\,(1 + \cos(2\pi X)),
\end{equation*}
with $A = 0.1$. Similar to the previous example, we consider two cases: (c) $N = 64, t_f = 0.1$ and (d) $N = 256, t_f = 0.1$, where the symbols have the usual meanings.

\begin{figure*}[tb]
    \centering
    \includegraphics[width = 0.95\textwidth]{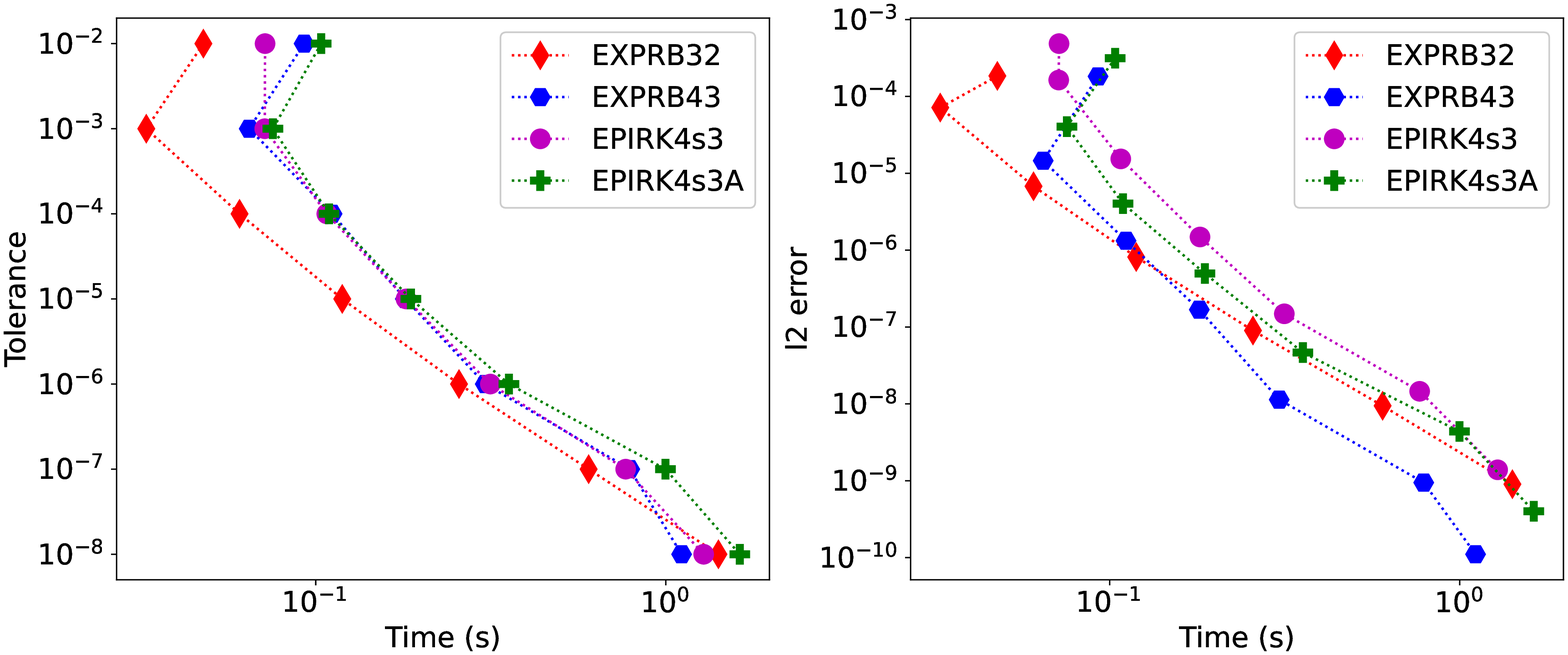} \\
    \includegraphics[width = 0.95\textwidth]{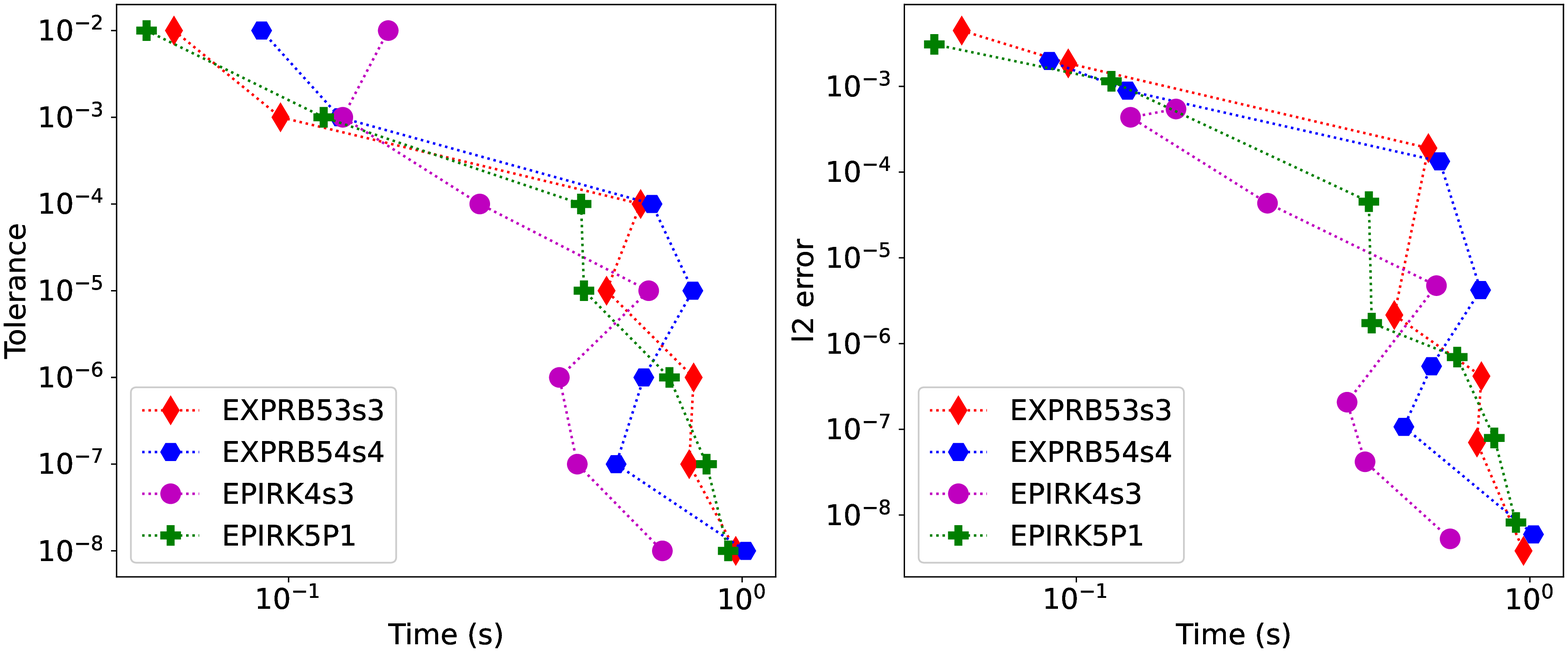}
    \caption{Illustration of the computational performance of a few embedded EXPRB and EPIRK integrators as a function of the user-defined tolerance (left panel) and the l2 norm of the global error incurred (right panel) for the Burgers' equation \eqref{eq:burgers}. The top and the bottom rows correspond to cases (a) and (b), respectively.}
    \label{fig:cost_bur}
\end{figure*}

\begin{figure*}[tb]
    \centering
    \includegraphics[width = 0.95\textwidth]{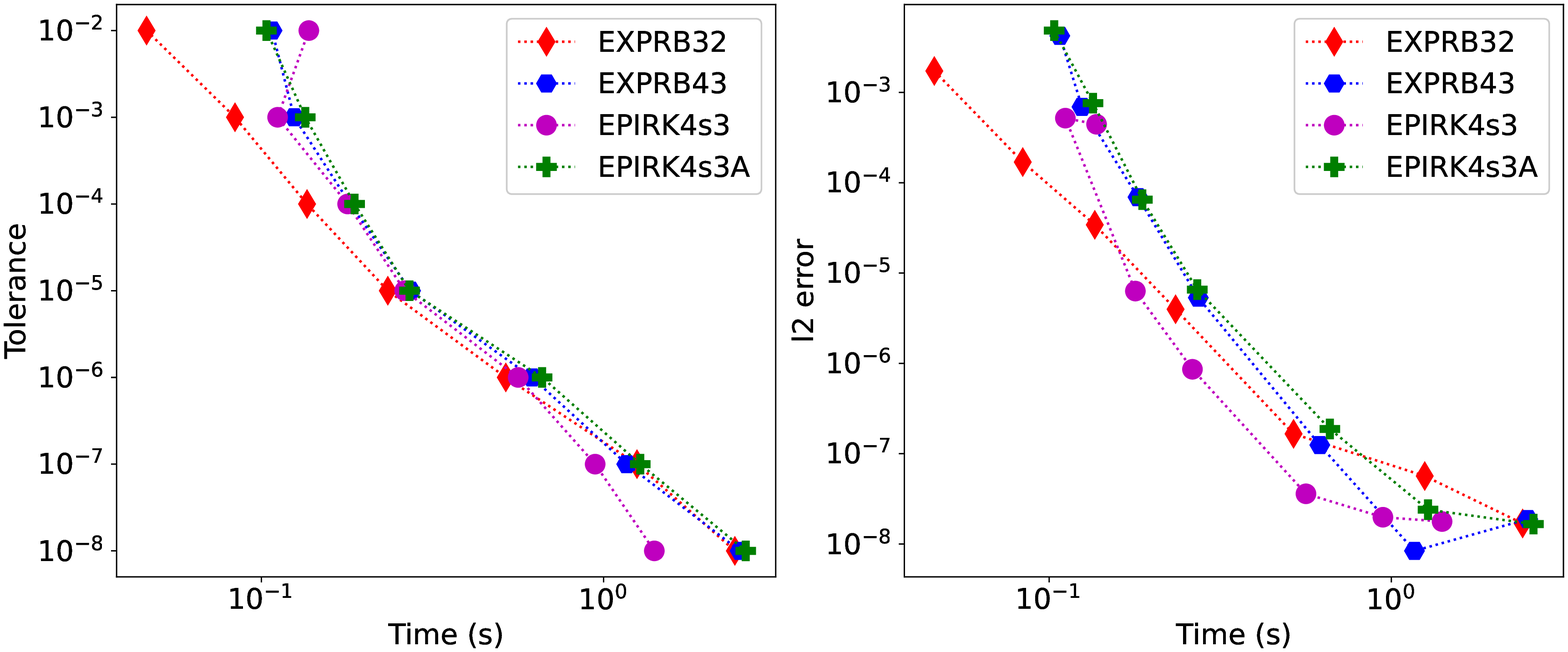} \\
    \includegraphics[width = 0.95\textwidth]{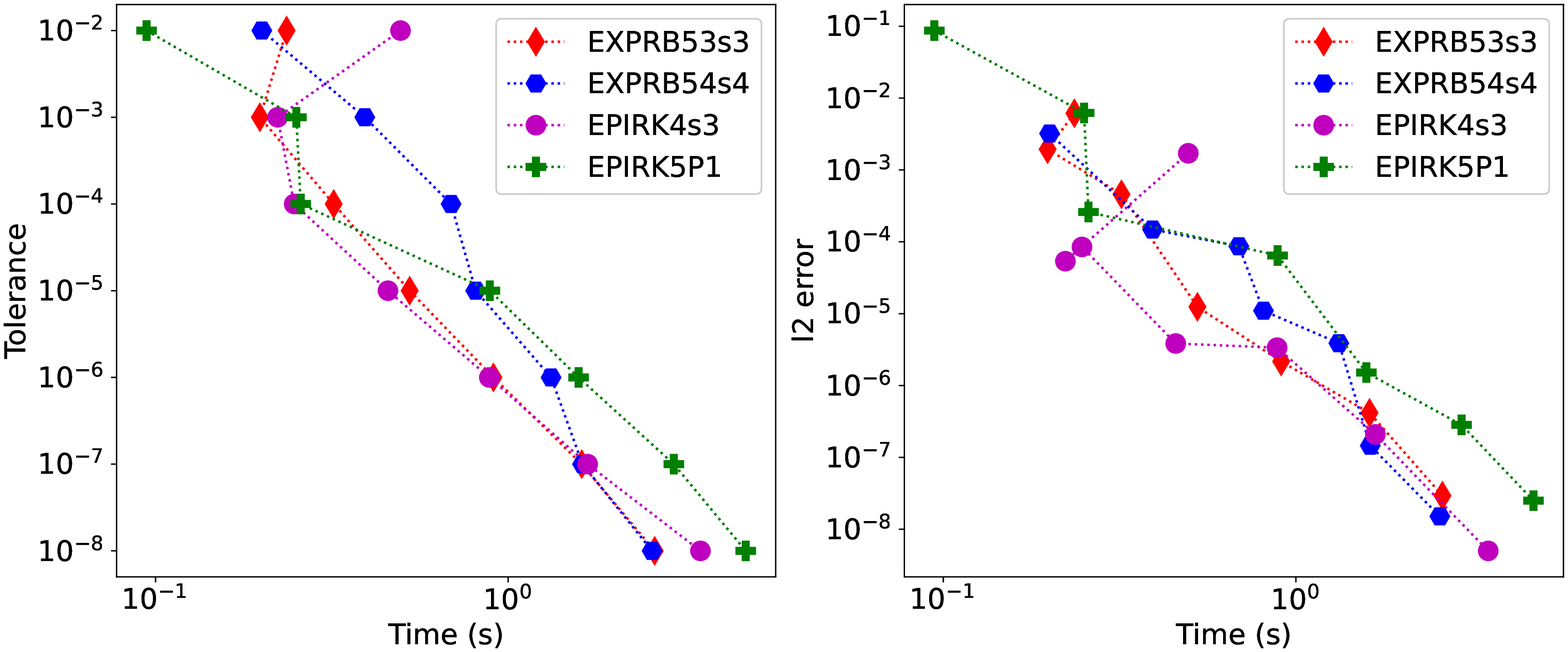}
    \caption{Illustration of the computational performance of a few embedded EXPRB and EPIRK integrators as a function of the user-defined tolerance (left panel) and the l2 norm of the global error incurred (right panel) for the Allen---Cahn equation \eqref{eq:ac}. The top and the bottom rows correspond to cases (c) and (d), respectively.}
    \label{fig:cost_ac}
\end{figure*}

We show the order of convergence for a selected number of integrators in Fig. \ref{fig:order} for the Burger's equation \eqref{eq:burgers} with $N = 64$ and $\eta = 200$. In Figs. \ref{fig:cost_bur} and \ref{fig:cost_ac}, we show the performance of a wide range of EXPRB and EPIRK integrators for a couple of representative problems (Eqs. \eqref{eq:burgers} and \eqref{eq:ac}) with variable step size implementation. Let us clearly state that, here, we do not investigate the performance of the different integrators. We simply demonstrate the performance of the different integrators available in \lexint rather than compare them.

%%%%%%%%%%%%%%%%%%%%%%%%%%%%%%%%%%%%%%%%%%%%%%%%%%%%%%%%%%%%%%%%%%%%%%%%%%%%%%%%%%

\section{Impact \& future aspects}

With the ever increasing need for high-resolution large-scale simulations in computational physics, there is a demand for ever more efficient and enhanced numerical algorithms. Efficiently integrating PDEs in time goes a long way in this regard. Exponential integrators have shown remarkable progress and promise in the last couple of decades. Various classes of exponential integrators have been shown to have superior performance to the traditional implicit and explicit methods for a wide range of problems usually considered in the mathematical literature. Additionally, their superiority have also been demonstrated for the MHD problems \cite{Tokman02, Einkemmer17, Deka22b}, kinetic plasma simulations \cite{Tuckmantel10, Dimarco11, Frenod15, Crouseilles18, Crouseilles20}, atmospheric and meteorological studies \cite{Clancy13, Gaudreault16, Mengalo18, Luan19, Schreiber19, Shashkin20, Brachet22, Pudykiewicz22}, and in different fields in engineering \cite{Rambeerich09, Michels14, Wang15, Tokman17a, Chen18, Chimmalgi19, Hammoud22}. 

Computing the exponential of a matrix constitutes a vital element for exponential integrators. Several approaches for the same have been outlined in the reviews by Moler and Van Loan \cite{Moler78, Moler03}. \texttt{EXPINT} \footnote{\url{http://www.math.ntnu.no/num/expint/matlab}} \cite{Berland07} is a \texttt{MATLAB} package for exponential integrators that computes the $\varphi_l(z)$ functions using a modified form of scaling and squaring (method 3 in \cite{Moler03}). The Krylov-subspace algorithm has become increasingly popular over the last few years owing to their ability to treat large system of matrices effectively. \texttt{EXPOKIT} \cite{Sidje98}, \texttt{phipm} \cite{Niesen12}, and \texttt{phipm\_simul\_iom} \footnote{\url{https://github.com/drreynolds/Phipm_simul_iom}} \cite{Luan19} are some the publicly available Krylov-based \texttt{MATLAB} software for efficiently computing the matrix exponential as well as computing the $\varphi_l(z)$ functions for exponential integrators. Further research in this field have shown that the method of polynomial interpolation \cite{Caliari04, Caliari07b} is highly competitive with, if not better than, the Krylov-based methods. \texttt{expleja} \footnote{\url{https://www.mathworks.com/matlabcentral/fileexchange/44039-matrix-exponential-times-a-vector}} is one of the first Leja-interpolation based \texttt{MATLAB} software that computes the matrix exponential times a vector or a matrix.

With the increasing popularity of exponential integrators in various fields in computational science, we provide an accessible framework for an efficient implementation of these methods. Whilst methods like Pad\'e approximation, squaring and scaling, or diagonalising the corresponding matrix and computing the exponential of the resulting eigenvalues work well for small matrices, these methods become prohibitive for large systems. Libraries based on (parts of) these methods are already available in \textsc{Python}. We provide a library based on the the Leja polynomial interpolation method that is highly favourable for the computation of the exponential-like functions of large systems of matrices. As part of the software package, we present a multitude of (Leja-based) exponential integrators (from the literature) for temporal integration of nonlinear PDEs with constant and variable step sizes. The \texttt{Matlab} version of the Leja interpolation method for exponential integrators has been appended to the Krylov-based \texttt{EPIC} library. Using \texttt{EPIC}, it has recently been shown by Gaudreault et al. \cite{Gaudreault18}, that an algorithm based on incomplete orthogonalisation of the basis vectors (\texttt{KIOPS}) may help in achieving a reasonable amount of improved performance over the state-of-the-art \texttt{phipm} algorithm. This publicly available package can be obtained from \url{https://faculty.ucmerced.edu/mtokman/#software}. We have used this package in our study of performance comparison of the Leja method with the \texttt{KIOPS} algorithm \cite{Deka22c} and to develop efficient integrators for anisotropic diffusion \cite{Deka22d}.

The release of the present package provides an effective implementation of the Leja-based method for people to get started on exploring such a method in a user-friendly environment. In the future, we will develop a parallel implementation of \lexint and include other exponential integrators that are designed specifically for parallel computing and for high-performance computing (such as GPUs). \lexint will be then implemented as a part of large software packages. As an example, in the near future, this package will be appended to the \texttt{PICARD} code \cite{Kissmann14} to solve the time-dependent cosmic-ray transport equation \cite{Strong07}.

%%%%%%%%%%%%%%%%%%%%%%%%%%%%%%%%%%%%%%%%%%%%%%%%%%%%%%%%%%%%%%%%%%%%%%%%%%%%%%%%%%

\section*{Acknowledgements}
This work is supported, in part, by the Austrian Science Fund (FWF) project id: P32143-N32. We would like to thank Marco Caliari for providing us with the code to compute Leja points.

%%%%%%%%%%%%%%%%%%%%%%%%%%%%%%%%%%%%%%%%%%%%%%%%%%%%%%%%%%%%%%%%%%%%%%%%%%%%%%%%%%

% \appendix

% \section{Appendix}
% \label{sec:appendix}

\bibliographystyle{elsarticle-num} 
\bibliography{ref}

\begin{thebibliography}{10}
\expandafter\ifx\csname url\endcsname\relax
  \def\url#1{\texttt{#1}}\fi
\expandafter\ifx\csname urlprefix\endcsname\relax\def\urlprefix{URL }\fi
\expandafter\ifx\csname href\endcsname\relax
  \def\href#1#2{#2} \def\path#1{#1}\fi

\bibitem{Ostermann10}
M.~Hochbruck, A.~Ostermann, {Exponential Integrators}, Acta Numer. 19 (2010)
  209 -- 286.
\newblock \href {https://doi.org/10.1017/S0962492910000048}
  {\path{doi:10.1017/S0962492910000048}}.

\bibitem{Deka22a}
P.~J. Deka, L.~Einkemmer, {Efficient adaptive step size control for exponential
  integrators}, Comput. Math. Appl. 123 (2022) 59--74.
\newblock \href {https://doi.org/10.1016/j.camwa.2022.07.011}
  {\path{doi:10.1016/j.camwa.2022.07.011}}.

\bibitem{Deka22b}
P.~J. Deka, L.~Einkemmer, {Exponential} {Integrators} for {Resistive}
  {Magnetohydrodynamics}: {Matrix-free} {Leja} {Interpolation} and {Efficient}
  {Adaptive} {Time} {Stepping}, ApJS 259~(2) (2022) 57.
\newblock \href {https://doi.org/10.3847/1538-4365/ac5177}
  {\path{doi:10.3847/1538-4365/ac5177}}.

\bibitem{Tokman16}
G.~{Rainwater}, M.~{Tokman}, {A new approach to constructing efficient stiffly
  accurate EPIRK methods}, J. Comput. Phys. 323 (2016) 283--309.
\newblock \href {https://doi.org/10.1016/j.jcp.2016.07.026}
  {\path{doi:10.1016/j.jcp.2016.07.026}}.

\bibitem{Caliari09}
M.~Caliari, A.~Ostermann, {Implementation of exponential Rosenbrock-type
  integrators}, Appl. Numer. Math. 59~(3) (2009) 568 -- 581.
\newblock \href {https://doi.org/10.1016/j.apnum.2008.03.021}
  {\path{doi:10.1016/j.apnum.2008.03.021}}.

\bibitem{Hochbruck09}
M.~Hochbruck, A.~Ostermann, J.~Schweitzer, {Exponential Rosenbrock-Type
  Methods}, SIAM J. Numer. Anal. 47~(1) (2009) 786--803.
\newblock \href {https://doi.org/10.1137/080717717}
  {\path{doi:10.1137/080717717}}.

\bibitem{Luan14}
V.~T. Luan, A.~Ostermann, {Exponential Rosenbrock methods of order five —
  construction, analysis and numerical comparisons}, J. Comput. Appl. Math. 255
  (2014) 417--431.
\newblock \href {https://doi.org/10.1016/j.cam.2013.04.041}
  {\path{doi:10.1016/j.cam.2013.04.041}}.

\bibitem{Tokman17a}
D.~L. Michels, V.~T. Luan, M.~Tokman, A stiffly accurate integrator for
  elastodynamic problems, ACM Trans. Graph. 36~(4) (2017).
\newblock \href {https://doi.org/10.1145/3072959.3073706}
  {\path{doi:10.1145/3072959.3073706}}.

\bibitem{Tokman17b}
G.~{Rainwater}, M.~{Tokman}, {Designing efficient exponential integrators with
  EPIRK framework}, in: International Conference of Numerical Analysis and
  Applied Mathematics (ICNAAM 2016), Vol. 1863 of American Institute of Physics
  Conference Series, 2017, p. 020007.
\newblock \href {https://doi.org/10.1063/1.4992153}
  {\path{doi:10.1063/1.4992153}}.

\bibitem{Tokman12}
M.~Tokman, J.~Loffeld, P.~Tranquilli, {New Adaptive Exponential Propagation
  Iterative Methods of Runge--Kutta Type}, SIAM J. Sci. Comput. 34~(5) (2012)
  A2650--A2669.
\newblock \href {https://doi.org/10.1137/110849961}
  {\path{doi:10.1137/110849961}}.

\bibitem{Pope63}
D.~A. Pope, An exponential method of numerical integration of ordinary
  differential equations, Commun. ACM 6~(8) (1963) 491–493.
\newblock \href {https://doi.org/10.1145/366707.367592}
  {\path{doi:10.1145/366707.367592}}.

\bibitem{Luan17}
V.~T. Luan, {Fourth-order two-stage explicit exponential integrators for
  time-dependent PDEs}, Appl. Numer. Math. 112 (2017) 91--103.
\newblock \href {https://doi.org/10.1016/j.apnum.2016.10.008}
  {\path{doi:10.1016/j.apnum.2016.10.008}}.

\bibitem{Sidje98}
R.~B. Sidje, Expokit: A software package for computing matrix exponentials, ACM
  Trans. Math. Softw. 24~(1) (1998) 130–156.
\newblock \href {https://doi.org/10.1145/285861.285868}
  {\path{doi:10.1145/285861.285868}}.

\bibitem{Moler03}
C.~Moler, C.~Van~Loan, {Nineteen Dubious Ways to Compute the Exponential of a
  Matrix, Twenty-Five Years Later}, SIAM Rev. 45~(1) (2003) 3--49.
\newblock \href {https://doi.org/10.1137/S00361445024180}
  {\path{doi:10.1137/S00361445024180}}.

\bibitem{Einkemmer17}
L.~Einkemmer, M.~Tokman, J.~Loffeld, On the performance of exponential
  integrators for problems in magnetohydrodynamics, J. Sci. Comput. 330 (2017)
  550--565.
\newblock \href {https://doi.org/10.1016/j.jcp.2016.11.027}
  {\path{doi:10.1016/j.jcp.2016.11.027}}.

\bibitem{Gaudreault18}
S.~{Gaudreault}, G.~{Rainwater}, M.~{Tokman}, {KIOPS: A fast adaptive Krylov
  subspace solver for exponential integrators}, J. Comput. Phys. 372 (2018)
  236--255.
\newblock \href {https://doi.org/10.1016/j.jcp.2018.06.026}
  {\path{doi:10.1016/j.jcp.2018.06.026}}.

\bibitem{Leja1957}
F.~Leja, \href{http://eudml.org/doc/208291}{Sur certaines suites liées aux
  ensembles plans et leur application à la représentation conforme}, Ann.
  Polon. Math. 4~(1) (1957) 8--13.
\newline\urlprefix\url{http://eudml.org/doc/208291}

\bibitem{Reichel90}
L.~Reichel, {Newton interpolation at Leja points}, BIT 30 (1990) 332 -- 346.
\newblock \href {https://doi.org/10.1007/BF02017352}
  {\path{doi:10.1007/BF02017352}}.

\bibitem{Baglama98}
J.~Baglama, D.~Calvetti, L.~Reichel, \href{http://eudml.org/doc/119747}{{Fast
  Leja Points}}, Electron. Trans. Numer. Anal. 7 (1998) 124 -- 140.
\newline\urlprefix\url{http://eudml.org/doc/119747}

\bibitem{Bergamaschi06}
L.~Bergamaschi, M.~Caliari, A.~Martinez, M.~Vianello, {Comparing Leja and
  Krylov Approximations of Large Scale Matrix Exponentials}, Proc. ICCS (2006)
  685--692\href {https://doi.org/10.1007/11758549_93}
  {\path{doi:10.1007/11758549_93}}.

\bibitem{Caliari07b}
M.~{Caliari}, M.~{Vianello}, L.~{Bergamaschi}, {The LEM exponential integrator
  for advection-diffusion-reaction equations}, J. Comput. Appl. Math. 210~(1-2)
  (2007) 56--63.
\newblock \href {https://doi.org/10.1016/j.cam.2006.10.055}
  {\path{doi:10.1016/j.cam.2006.10.055}}.

\bibitem{Caliari04}
M.~Caliari, M.~Vianello, L.~Bergamaschi, {Interpolating discrete
  advection–diffusion propagators at Leja sequences}, J. Comput. Appl. Math.
  172~(1) (2004) 79 -- 99.
\newblock \href {https://doi.org/10.1016/j.cam.2003.11.015}
  {\path{doi:10.1016/j.cam.2003.11.015}}.

\bibitem{Caliari14}
M.~Caliari, P.~Kandolf, A.~Ostermann, S.~Rainer, Comparison of software for
  computing the action of the matrix exponential, BIT Numer. Math. 54 (2014)
  113 -- 128.
\newblock \href {https://doi.org/10.1007/s10543-013-0446-0}
  {\path{doi:10.1007/s10543-013-0446-0}}.

\bibitem{Einkemmer18}
L.~Einkemmer, An adaptive step size controller for iterative implicit methods,
  Appl. Numer. Math. 132 (2018) 182 -- 204.
\newblock \href {https://doi.org/10.1016/j.apnum.2018.06.002}
  {\path{doi:10.1016/j.apnum.2018.06.002}}.

\bibitem{Tokman02}
M.~Tokman, P.~M. Bellan, {Three-dimensional Model of the Structure and
  Evolution of Coronal Mass Ejections}, ApJ 567~(2) (2002) 1202--1210.
\newblock \href {https://doi.org/10.1086/338699} {\path{doi:10.1086/338699}}.

\bibitem{Tuckmantel10}
T.~Tuckmantel, A.~Pukhov, J.~Liljo, M.~Hochbruck, {Three-Dimensional
  Relativistic Particle-in-Cell Hybrid Code Based on an Exponential
  Integrator}, IEEE Transactions on Plasma Science 38~(9) (2010) 2383--2389.
\newblock \href {https://doi.org/10.1109/TPS.2010.2056706}
  {\path{doi:10.1109/TPS.2010.2056706}}.

\bibitem{Dimarco11}
G.~Dimarco, L.~Pareschi, {Exponential Runge–Kutta Methods for Stiff Kinetic
  Equations}, SIAM Journal on Numerical Analysis 49~(5) (2011) 2057--2077.
\newblock \href {https://doi.org/10.1137/100811052}
  {\path{doi:10.1137/100811052}}.

\bibitem{Frenod15}
E.~Frenod, S.~A. Hirstoaga, M.~Lutz, E.~Sonnendrücker, {Long Time Behaviour of
  an Exponential Integrator for a Vlasov-Poisson System with Strong Magnetic
  Field}, CiCP 18~(2) (2015) 263–296.
\newblock \href {https://doi.org/10.4208/cicp.070214.160115a}
  {\path{doi:10.4208/cicp.070214.160115a}}.

\bibitem{Crouseilles18}
N.~Crouseilles, L.~Einkemmer, M.~Prugger, An exponential integrator for the
  drift-kinetic model, Comput. Phys. Commun. 224 (2018) 144--153.
\newblock \href {https://doi.org/10.1016/j.cpc.2017.11.003}
  {\path{doi:10.1016/j.cpc.2017.11.003}}.

\bibitem{Crouseilles20}
N.~Crouseilles, L.~Einkemmer, J.~Massot, Exponential methods for solving
  hyperbolic problems with application to collisionless kinetic equations, J.
  Comput. Phys. 420 (2020) 109688.
\newblock \href {https://doi.org/10.1016/j.jcp.2020.109688}
  {\path{doi:10.1016/j.jcp.2020.109688}}.

\bibitem{Clancy13}
C.~Clancy, J.~A. Pudykiewicz, On the use of exponential time integration
  methods in atmospheric models, Tellus A: Dynamic Meteorology and Oceanography
  65~(1) (2013) 20898.
\newblock \href {https://doi.org/10.3402/tellusa.v65i0.20898}
  {\path{doi:10.3402/tellusa.v65i0.20898}}.

\bibitem{Gaudreault16}
S.~Gaudreault, J.~A. Pudykiewicz, An efficient exponential time integration
  method for the numerical solution of the shallow water equations on the
  sphere, J. Comput. Phys. 322 (2016) 827--848.
\newblock \href {https://doi.org/10.1016/j.jcp.2016.07.012}
  {\path{doi:10.1016/j.jcp.2016.07.012}}.

\bibitem{Mengalo18}
G.~Mengaldo, A.~Wyszogrodzki, M.~Diamantakis, S.-J. Lock, F.~X. Giraldo, N.~P.
  Wedi, {Current and Emerging Time-Integration Strategies in Global Numerical
  Weather and Climate Prediction}, Arch. Computat. Methods Eng. 26 (2019)
  663–684.
\newblock \href {https://doi.org/doi.org/10.1007/s11831-018-9261-8}
  {\path{doi:doi.org/10.1007/s11831-018-9261-8}}.

\bibitem{Luan19}
V.~T. Luan, J.~A. Pudykiewicz, D.~R. Reynolds, Further development of efficient
  and accurate time integration schemes for meteorological models, J. Comput.
  Phys. 376 (2019) 817--837.
\newblock \href {https://doi.org/10.1016/j.jcp.2018.10.018}
  {\path{doi:10.1016/j.jcp.2018.10.018}}.

\bibitem{Schreiber19}
M.~Schreiber, N.~Schaeffer, R.~Loft, Exponential integrators with
  parallel-in-time rational approximations for the shallow-water equations on
  the rotating sphere, Parallel Computing 85 (2019) 56--65.
\newblock \href {https://doi.org/10.1016/j.parco.2019.01.005}
  {\path{doi:10.1016/j.parco.2019.01.005}}.

\bibitem{Shashkin20}
V.~Shashkin, G.~Goyman, {Parallel Efficiency of Time-Integration Strategies for
  the Next Generation Global Weather Prediction Model}, in: V.~Voevodin,
  S.~Sobolev (Eds.), Supercomputing, Springer International Publishing, Cham,
  2020, pp. 285--296.

\bibitem{Brachet22}
M.~Brachet, L.~Debreu, C.~Eldred, Comparison of exponential integrators and
  traditional time integration schemes for the shallow water equations, Appl.
  Numer. Math. 180 (2022) 55--84.
\newblock \href {https://doi.org/10.1016/j.apnum.2022.05.006}
  {\path{doi:10.1016/j.apnum.2022.05.006}}.

\bibitem{Pudykiewicz22}
J.~A. Pudykiewicz, C.~Clancy, Convection experiments with the exponential time
  integration scheme, J. Comput. Phys. 449 (2022) 110803.
\newblock \href {https://doi.org/10.1016/j.jcp.2021.110803}
  {\path{doi:10.1016/j.jcp.2021.110803}}.

\bibitem{Rambeerich09}
N.~{Rambeerich}, D.~Y. {Tangman}, A.~{Gopaul}, M.~{Bhuruth}, {Exponential time
  integration for fast finite element solutions of some financial engineering
  problems}, J. Comput. Appl. Math. 224~(2) (2009) 668--678.
\newblock \href {https://doi.org/10.1016/j.cam.2008.05.047}
  {\path{doi:10.1016/j.cam.2008.05.047}}.

\bibitem{Michels14}
D.~L. Michels, G.~A. Sobottka, A.~G. Weber, {Exponential Integrators for Stiff
  Elastodynamic Problems}, ACM Trans. Graph. 33~(1) (2014).
\newblock \href {https://doi.org/10.1145/2508462} {\path{doi:10.1145/2508462}}.

\bibitem{Wang15}
C.~Wang, X.~Fu, P.~Li, J.~Wu, Accurate dense output formula for exponential
  integrators using the scaling and squaring method, Appl. Math. Lett. 43
  (2015) 101--107.
\newblock \href {https://doi.org/10.1016/j.aml.2014.12.008}
  {\path{doi:10.1016/j.aml.2014.12.008}}.

\bibitem{Chen18}
Y.~J. Chen, U.~M. Ascher, D.~K. Pai, {Exponential Rosenbrock-Euler Integrators
  for Elastodynamic Simulation}, Vol. 24(10), IEEE, 2018.
\newblock \href {https://doi.org/10.1109/TVCG.2017.2768532}
  {\path{doi:10.1109/TVCG.2017.2768532}}.

\bibitem{Chimmalgi19}
S.~Chimmalgi, P.~J. Prins, S.~Wahls, {Fast Nonlinear Fourier Transform
  Algorithms Using Higher Order Exponential Integrators}, IEEE Access 7 (2019)
  145161--145176.
\newblock \href {https://doi.org/10.1109/ACCESS.2019.2945480}
  {\path{doi:10.1109/ACCESS.2019.2945480}}.

\bibitem{Hammoud22}
B.~{Hammoud}, L.~{Olivieri}, L.~{Righetti}, J.~{Carpentier}, A.~{Del Prete},
  {Exponential integration for efficient and accurate multibody simulation with
  stiff viscoelastic contacts}, Multibody Syst. Dyn. 54 (2022) 443--460.
\newblock \href {https://doi.org/10.1007/s11044-022-09818-z}
  {\path{doi:10.1007/s11044-022-09818-z}}.

\bibitem{Moler78}
C.~B. Moler, C.~V. Loan, {Nineteen Dubious Ways to Compute the Exponential of a
  Matrix}, Siam Rev. 20 (1978) 801--836.
\newblock \href {https://doi.org/10.1137/1020098} {\path{doi:10.1137/1020098}}.

\bibitem{Berland07}
H.~Berland, B.~Skaflestad, W.~M. Wright, {EXPINT---A MATLAB Package for
  Exponential Integrators}, ACM Trans. Math. Softw. 33~(1) (2007) 4–es.
\newblock \href {https://doi.org/10.1145/1206040.1206044}
  {\path{doi:10.1145/1206040.1206044}}.

\bibitem{Niesen12}
J.~Niesen, W.~M. Wright, {Algorithm 919: A Krylov Subspace Algorithm for
  Evaluating the $\varphi$-Functions Appearing in Exponential Integrators}, ACM
  Trans. Math. Softw. 38~(3) (2012).
\newblock \href {https://doi.org/10.1145/2168773.2168781}
  {\path{doi:10.1145/2168773.2168781}}.

\bibitem{Deka22c}
P.~J. {Deka}, M.~{Tokman}, L.~{Einkemmer}, {A comparison of Leja- and
  Krylov-based iterative schemes for Exponential Integrators}, arXiv e-prints
  (2022) arXiv:2211.08948\href {http://arxiv.org/abs/2211.08948}
  {\path{arXiv:2211.08948}}.

\bibitem{Deka22d}
P.~J. {Deka}, L.~{Einkemmer}, R.~{Kissmann}, {Exponential methods for
  anisotropic diffusion}, arXiv e-prints (2022) arXiv:2211.08953\href
  {http://arxiv.org/abs/2211.08953} {\path{arXiv:2211.08953}}.

\bibitem{Kissmann14}
R.~Kissmann, {\textsc{PICARD}: A novel code for the Galactic Cosmic Ray
  propagation problem}, Astropart. Phys. 55 (2014) 37--50.
\newblock \href {https://doi.org/10.1016/j.astropartphys.2014.02.002}
  {\path{doi:10.1016/j.astropartphys.2014.02.002}}.

\bibitem{Strong07}
A.~W. Strong, I.~V. Moskalenko, V.~S. Ptuskin, {Cosmic-Ray Propagation and
  Interactions in the Galaxy}, Annu. Rev. Nucl. Part. Sci. 57~(1) (2007)
  285--327.
\newblock \href {https://doi.org/10.1146/annurev.nucl.57.090506.123011}
  {\path{doi:10.1146/annurev.nucl.57.090506.123011}}.

\end{thebibliography}





% \end{thebibliography}
\end{document}